\documentclass{amsart}
\usepackage{mathtools, microtype, mathrsfs, xfrac, hyperref, amssymb, enumerate, xspace}
\usepackage[utf8]{inputenc}						
\usepackage{tikz-cd}
\usepackage{xcolor}

\usepackage[style=alphabetic, backend=biber]{biblatex}
\bibliography{real}

\numberwithin{equation}{section}

\numberwithin{subsection}{section}

\allowdisplaybreaks[1]

\newenvironment{enumeratei}
{\begin{enumerate}[\upshape (i)]}
{\end{enumerate}}

\newenvironment{enumerate1}
{\begin{enumerate}[\upshape (1)]}
{\end{enumerate}}

\newtheorem*{namedtheorem}{\theoremname}
\newcommand{\theoremname}{testing}

\newtheorem{proposition-definition}[theorem]
{Proposition-Definition}

\newtheorem*{theorem*}{Theorem}

\theoremstyle{definition}

\theoremstyle{remark}

\renewcommand{\mathcal}{\mathscr}

\renewcommand\AA{\mathbb{A}} 
\newcommand\CC{\mathbb{C}} 
 
\newcommand\GG{\mathbb{G}}

 \newcommand\RR{\mathbb{R}}

 \newcommand\ZZ{\mathbb{Z}}

\newcommand\arr{\ifinner\to\else\longrightarrow\fi}

\newcommand\arrto{\ifinner\mapsto\else\longmapsto\fi}

\renewcommand\H{\operatorname{H}}

\newcommand{\eqdef}{\mathrel{\smash{\overset{\mathrm{\scriptscriptstyle def}} =}}}

\def\displaytimes_#1{\mathrel{\mathop{\times}\limits_{#1}}}

\def\displayotimes_#1{\mathrel{\mathop{\bigotimes}\limits_{#1}}}

\renewcommand\projlim{\varprojlim}

\newlength{\ignora}

\renewcommand{\setminus}{\smallsetminus}

\renewcommand\projlim{\varprojlim}

\DeclareFontFamily{U}{mathx}{\hyphenchar\font45}
\DeclareFontShape{U}{mathx}{m}{n}{
      <5> <6> <7> <8> <9> <10>
      <10.95> <12> <14.4> <17.28> <20.74> <24.88>
      mathx10
      }{}
\DeclareSymbolFont{mathx}{U}{mathx}{m}{n}
\DeclareFontSubstitution{U}{mathx}{m}{n}
\DeclareMathAccent{\widecheck}{0}{mathx}{"71}
\DeclareMathAccent{\wideparen}{0}{mathx}{"75}

\renewcommand{\epsilon}{\varepsilon}

\newcommand{\sect}{\operatorname{Sec}}

\begin{document}

\title[The real section conjecture]{An elementary approach to Stix's proof of the real section conjecture}

\author[Vistoli]{Angelo Vistoli}
\address{Scuola Normale Superiore\\Piazza dei Cavalieri 7\\
56126 Pisa\\ Italy}
\email[vistoli]{angelo.vistoli@sns.it}

\author[Bresciani]{Giulio Bresciani}
\address{Freie Universit\"at Berlin, Arnimallee 3, 14195 Berlin, Germany}
\email[bresciani]{gbresciani@zedat.fu-berlin.de}

\thanks{The second author is supported by the DFG Priority Program "Homotopy Theory and Algebraic Geometry" SPP 1786}


\maketitle


\section{Introduction}
Given a geometrically connected variety $X$ over a field $k$ with separable closure $\bar{k}$ and a geometric point $\overline{x} \in X(\overline{k})$, there is a short exact sequence of \'etale fundamental groups
\[1 \to \pi_{1}(X_{\bar{k}},\overline{x}) \to \pi_{1}(X,\overline{x}) \to \operatorname{Gal}(\bar{k}/k) \to 1.\]
A rational point $x\in X(k)$ yields by functoriality a section of this short exact sequence which is well defined up to conjugacy by elements of $\pi_{1}(X_{\bar{k}})$. If we call $\sect(X/k)$ the space of sections modulo the action of conjugation by $\pi_{1}(X_{\bar{k}})$, this gives a natural map
\[X(k)\to \sect(X/k)\]
called the \emph{section map}. Grothendieck's famous section conjecture predicts that the section map is a bijection if $k$ is a number field and $X$ is a smooth, proper curve of genus at least 2.

If $k=\RR$, the section map factors through a map
\[\pi_{0}\bigl(X(\RR)\bigr) \to \sect(X/k)\]
where $\pi_{0}\bigl(X(\RR)\bigr)$ is the set of connected components of $X(\RR)$, see \cite[A.4]{st10}.

When this map is a bijection, we say that \emph{the real section conjecture holds for $X$}. S.~Mochizuki proved in \cite{mo03} that the real section conjecture holds for smooth, quasiprojective hyperbolic curves, and for abelian varieties, using a theorem of Cox which computes the \'etale $\ZZ/2$-cohomology of $X$ in terms of the singular $\ZZ/2$-cohomology of $X(\RR)$. After Mochizuki, several different proofs were found, at least for curves: J.~Stix \cite{st10} used a theorem of Witt which computes the Brauer group of the curve in terms of the singular cohomology of $X(\RR)$, A.~P\'al \cite{pal11} used topological methods in the spirit of Smith’s fixed-point theorem, K.~Wickelgren \cite{wick14} used an obstruction coming from the maximal 2-nilpotent quotient of the fundamental group. For curves of genus $g\ge 1$ without real points, H.~Esnault and O.~Wittenberg gave a proof using only Tsen's theorem, see \cite[Remark 3.7(iv)]{ew09}.

The purpose of this note is to give an elementary proof of Mochizuki's result, using only Kummer theory, basic topological properties of real curves (a modern account of these properties is given by \cite{gh81}), and a well-known argument which allows to reduce to the case of curves without real points using Tamagawa's neighborhoods of a section. Our argument follows the lines of Stix's: still, we think it is worth observing that no kind of advanced technology is needed. 

\begin{theorem*}
	The real section conjecture holds for the following two classes of geometrically connected real varieties.
\begin{enumerate1}

\item Torsors for semiabelian varieties.

\item Smooth quasiprojective curves, with the exception of the projective curves of genus $0$ without real points.

\end{enumerate1}
\end{theorem*}

It should be remarked that, although Mochizuki only states his result for abelian varieties, in fact his proof goes through for semiabelian varieties.

\section{The proof}
The proof is in three steps. First we prove the case of torsors for semiabelian varieties; the argument is basically Kummer theory.

The second case is that of projective curves without real points: this follows from the preceding case.

Finally, in the case of quasiprojective curves, surjectivity is proved from the second case, using Tamagawa's neighborhoods of a section, while for injectivity we use the first case and Rosenlicht's generalized jacobians.

Before starting, let us recall briefly Tamagawa's idea of neighborhoods of a section. A finite \'etale cover $Y\to X$ is a neighborhood of a section $s\in \sect(X/k)$ if $s$ lifts to $\sect(Y/k)$. Neighborhoods of $s$ form a projective system $\{X_{i}\}$. We will use the following two facts
\begin{enumeratei}
	\item a section is geometric (i.e. in the image of the section map) if and only if the projective limit $\projlim_{i}X_{i}$ has a rational point,
	\item the base change of $\{X_{i}\}$ to $\overline{k}$ is cofinal in the system of all finite \'etale covers of $X_{\overline{k}}$.
\end{enumeratei}
See \cite[Chapter 4]{st13} for details.

\subsection*{Torsors for semiabelian varieties}

If $A$ is a semiabelian variety over $\RR$, call $TA = \projlim_{n}A(\CC)[n]$ the Tate module of $A$. We have that $TA = \pi_{1}(A_{\CC})$ is a $\operatorname{Gal}(\CC/\RR)$-module and $\sect(A/\RR)=\H^{1}(\RR,TA)$, see for instance \cite[Corollary 71]{st13}. Kummer theory gives us a map
\[\delta\colon A(\RR)\to\H^{1}(\RR,TA)=\sect(A/\RR)\]
which coincides with the section map.

Let us prove injectivity. We have that $A(\RR)$ is a real, abelian Lie group, thus it is an extension of a finite group $F=\pi_{0}(A(\RR))$ by a connected group $C$. We already know that $\delta$ factors through $\pi_{0}(A(\RR))=F$, and Kummer theory tells us that the kernel of $\delta$ consists of divisible elements. Since $F$ is finite, it has no non-trivial divisible elements and thus $\pi_{0}(A(\RR))=F\to \H^{1}(\RR,TA)$ is injective.

Let us prove surjectivity. For every $n$, call $\phi_{n}$ the composition
\[\phi_{n}\colon \H^{1}(\RR,TA)\to\H^{1}(\RR,A[n])\to\H^{1}(\RR,A).\]
Kummer theory tells us that an element $a\in\H^{1}(\RR,TA)$ is in the image of $\delta$ if and only if $\phi_{n}(a)=0$ for every $n$. Recall now that $\H^{1}(\RR,A)$ is $2$-torsion, since $|\operatorname{Gal}(\CC/\RR)|=2$. Since $\phi_{n}(a)=2\phi_{2n}(a)=0$ for every $n$, we have that $\delta$ is surjective. This concludes the proof for semiabelian varieties.

If $B$ is a non-trivial torsor for a semiabelian variety $A$, we want to prove that $\sect(B/\RR)=\emptyset$. Suppose by contradiction that a section $s\in\sect(B/\RR)$ exists. Since $\H^{1}(\RR,A)$ is $2$-torsion, $B/A[2]\simeq A$, call $p\colon B\to B/A[2]=A$ the projection. By the preceding case, $p_{*}s\in\sect(A/\RR)$ is geometric, and $B$ is a neighborhood of $p_{*}s$. But this is a contradiction, since neighborhoods of geometric sections always have rational points by property (i) of neighborhoods.

\subsection*{Projective curves without real points}

Let $X$ be a smooth, projective and geometrically connected curve over $\RR$ without real points. By hypothesis, $X$ has positive genus. Suppose by contradiction that a section $s\in\sect(X/\RR)$ exists.

By property (ii), there exists a neighborhood $Y\to X$ of $s$ which is a finite \'etale cover of even degree, thus $Y$ has odd genus and no real points. Since $Y$ has odd genus and no real points, then $\operatorname{\underline{Pic}^{1}_{Y}}$ is a non-trivial abelian torsor, see \cite[Proposition 3.3.2]{gh81}. This implies that $\sect(\operatorname{\underline{Pic}^{1}_{Y}}/\RR)$ is empty by the preceding case. Since there is a natural map $\sect(Y/\RR)\to\sect(\operatorname{\underline{Pic}^{1}_{Y}}/\RR)$, then $\sect(Y/\RR)$ is empty, too. This is a contradiction, since $s$ lifts to $\sect(Y/\RR)$.

\subsection*{Quasiprojective curves}

Let $X$ be any smooth geometrically connected curve over $\RR$; let $\overline{X}$ be a smooth projective curve containing $X$ as a dense open subscheme, and write $X = \overline{X} \setminus D$, where $D$ is an effective reduced divisor on $\overline{X}$ of positive degree. 

Assume that the genus of $\overline{X}$ is $0$. If $\deg D = 1$, then $X = \AA^{1}_{\RR}$, and the result is obvious. If $\deg D = 2$, then it is easily seen that $X$ is a torsor under a $1$-dimensional torus, and we are back to the preceding case. So we may assume that either the genus of $\overline{X}$ is at least $1$, or $\deg D \geq 3$.

\subsubsection*{Injectivity of the section map} 

In order to prove injectivity, we may assume that $X(\RR)$ is non-empty, fix an origin $x\in X(\RR)$. Let $X\to J$ be Rosenlicht's generalized Jacobian, i.e. the semiabelian Albanese variety (see \cite[Ch.~5]{serre-algebraic-groups-class-fields}). Since we already know that the section map is injective on connected components for semiabelian varieties and in general the section map $Y(k)\to\sect(Y/k)$ is a natural transformation of functors, it is enough to prove that $\pi_{0}\bigl(X(\RR)\bigr)\to\pi_{0}(J(\RR))$ is injective.

Let $\overline{X}$ be the smooth compactification of $X$, $D=\overline{X}\setminus X$ and $\overline{X}\to\bar{J}$ the Jacobian, thanks to \cite[Proposition 4.2]{gh81} we have that $\pi_{0}(\overline{X}(\RR))\to\pi_{0}(\bar{J}(\RR))$ is injective. Since we know that injectivity holds in the projective case, to prove the general case we may focus on a single connected component $U$ of $\overline{X}(\RR)$, and suppose that $D\subset U$. We have that $U$ is homeomorphic to $S^{1}$, since this is the only compact, connected $1$-dimensional manifold. If $\deg D\le 1$, then $U\setminus D$ is connected, thus there is nothing to prove. The case $\deg D\ge 3$ follows from the case $\deg D=2$: if $C_{1},C_{2}\subset U\setminus D=S^{1}\setminus D$ are two different connected components, we may choose two points $p,q\in D$ such that $C_{1},C_{2}$ are contained in different connected components of $S^{1}\setminus\{p,q\}$, thus they map to different connected components in the Jacobian of $\overline{X}\setminus\{p,q\}$ thanks to the $\deg(D)=2$ case. Let us thus assume that $D=\{p,q\}\subset U$.

If $D=\{p,q\}$, then $J$ is an extension of $\bar{J}$ by $\GG_{m}$, and this extension is associated with the line bundle $L=[p-q]\in\hat{\bar{J}}(\RR)=\bar{J}(\RR)$. We have that $J\to\bar{J}$ is thus the $\GG_{m}$-torsor associated with $L$, i.e. $J$ is $L$ minus the $0$-section $J\to L$. Let $V\subset \bar{J}(\RR)$ be the connected component containing $U\subset \overline{X}(\RR)\subset \bar{J}(\RR)$. The restriction of $J(\RR)\to \bar{J}(\RR)$ to $V$ is thus an $\RR^{*}$-bundle $\tilde{V}\to V$. We want to show that the two connected components of $U\setminus \{p,q\}$ map to different connected components of $\tilde{V}$, and this in turn is equivalent to showing that $\tilde{V}$ is disconnected.

Now, recall that $L$ is associated with $[p-q]\in\bar{J}(\RR)$. Since $p,q\in\overline{X}(\RR)$ are in the same connected component, it follows that $[p-q]\in\bar{J}(\RR)$ is in the connected component of the identity. This implies that the $\RR^{*}$-bundle $\tilde{V}\to V$ is topologically trivial and thus disconnected.
 
\subsubsection*{Surjectivity}

Fix a section $s \in \sect(X/\RR)$, and consider the corresponding tower of neighborhoods $\{X_{i}\}$. It is enough to prove that $\projlim X_{i}(\RR) \neq \emptyset$. Notice that, because of the conditions on $X$, $D$ above, we may assume that the genus of $\overline{X}_{i}$ is positive, even if the genus of $\overline{X}$ is $0$.

For each $i$ call $\overline{X}_{i}$ a smooth projective compactification of $X_{i}$; the tower $\{X_{i}\}$ extends to a tower $\{\overline{X_{i}}\}$ mapping to $\overline{X}$. By construction, $\sect(X_{i}/\RR) \neq \emptyset$ for all $i$, hence $\sect(\overline{X}_{i}/\RR) \neq \emptyset$ for all $i$. Then, because of the previous case, we have $\overline{X_{i}}(\RR) \neq \emptyset$; since the $\overline{X_{i}}(\RR)$ are compact, this implies that $\projlim \overline{X}_{i}(\RR) \neq \emptyset$. Take a point $\pi$ in $\projlim \overline{X}_{i}(\RR)$; call $\overline{C}_{i}$ the connected component of $X_{i}(\RR)$ containing the image of $\pi$, and, similarly, call $\overline{C}$ the connected component of $\overline{X}(\RR)$ containing the image of $\pi$. Set $C_{i} \eqdef \overline{C}_{i} \cap X_{i}(\RR)$ and $C \eqdef \overline{C} \cap X(\RR)$. Clearly $C_{i}$ is the inverse image of $C$ in $\overline{C}_{i}$. Then the $C_{i}$ form a projective system mapping to $C$; the maps $C_{i} \arr C$ are surjective with finite fibers. If $p \in C$, the fiber over $p$ form a projective system of nonempty finite sets, which has nonempty projective limit, showing that $\projlim X_{i}(\RR) \neq \emptyset$, as claimed.

\printbibliography

\end{document}